\documentclass[amsppt,11pt]{amsart}
\usepackage{amssymb}
\usepackage{amsmath}

\newcommand{\Label}[1]{\label{#1}\hspace{.3cm}\fbox{\rm #1}\hspace{.3cm}}
\renewcommand{\Label}{\label} %% COMMENT THIS LINE OUT TO INCLUDE EDDY LABELS
\newtheorem {theorem}{Theorem}[section]
\newtheorem {proposition}[theorem]{Proposition}
\newtheorem {lemma}[theorem]{Lemma}
\newtheorem {corollary}[theorem]{Corollary}

\newcommand{\eqnsection}{
   \renewcommand{\theequation}{\thesection.\arabic{equation}}
   \makeatletter
   \csname @addtoreset\endcsname{equation}{section}
   \makeatother}
\def \med{ { \mbox{\rm  Med}}}

\newcommand{\be}{{\begin{equation}}}
\newcommand{\ee}{{\end{equation}}}
\def \bt{\begin{theorem}}
\def \et{\end{theorem}}
\def \bea{\begin{eqnarray}}
\def \eea{\end{eqnarray}}
\def \bas{\begin{eqnarray*}}
\def \eas{\end{eqnarray*}}

\def \ep{\eps}

\newcommand{\eps}{\varepsilon}

\newcommand{\ls}[1]
   {\dimen0=\fontdimen6\the\font \lineskip=#1\dimen0
\advance\lineskip.5\fontdimen5\the\font \advance\lineskip-\dimen0
\lineskiplimit=.9\lineskip \baselineskip=\lineskip
\advance\baselineskip\dimen0 \normallineskip\lineskip
\normallineskiplimit\lineskiplimit \normalbaselineskip\baselineskip
\ignorespaces }

\def \R{{\Bbb{R}}}
\def \RRR{{\bf R}}
\def \G{{\bf G}}
\def \Z{{\Bbb{Z}}}

\def \AA{{\mathcal A}}
\def \BB{{\mathcal B}}
\def \CC{{\mathcal C}}
\def \DD{{\mathcal D}}

\def \EE{{\mathcal E}}
\def \FF{{\mathcal F}}
\def \GG{{\mathcal G}}
\def \HH{{\mathcal H}}

\def \VV{{\mathcal V}}

\def \({\left(}
\def \){\right)}

\def \bc{\begin{center} }
\def \ec{\end{center} }
\def\Bbb{\mathbb}

\begin{document}

\eqnsection
\newcommand{\Ini}{{I_{n,i}}}
\newcommand{\reals}{{\Bbb{R}}}
\newcommand{\gff}{{\mathcal X}}
\newcommand{\tgff}{{\mathcal Y}}
\newcommand{\brw}{{\mathcal R}}
\newcommand{\mbrw}{{\mathcal S}}
\newcommand{\F}{{\mathcal F}}
\newcommand{\D}{{\mathcal D}}
\newcommand{\Fn}{{{\mathcal F}_n}}
\newcommand{\Gn}{{{\mathcal G}_n}}
\newcommand{\Hn}{{{\mathcal H}_n}}
\newcommand{\Fp}{{{\mathcal F}^p}}
\newcommand{\Gp}{{{\mathcal G}^p}}
\newcommand{\PPP}{{\mathbf P}}
\newcommand{\Pop}{{P\otimes \PPP}}
\newcommand{\hm}{\HH^\varphi}
\newcommand{\nuw}{{\nu^W}}
\newcommand{\ths}{{\theta^*}}
\newcommand{\beq}[1]{\begin{equation}\label{#1}}
\newcommand{\eeq}{\end{equation}}
\newcommand{\integers}{{\rm I\!N}}
\newcommand{\DT}{{D_{\Bbb T^2}}}
\newcommand{\pDT}{{\partial \DT}}
\newcommand{\E}{{\Bbb E}}
\newcommand{\te}{{\tilde{\delta}}}
\newcommand{\tI}{{\tilde{I}}}
\newcommand{\epn}{{\ep_n}}
\def\var{{\rm Var}}
\def\cov{{\rm Cov}}
\def\one{{\bf 1}}
\def\leb{{\mathcal L}eb}
\def\Ho{{\mbox{\sf H\"older}}}  %% This is H\"older exponent in formulas.
\def\thi{{\mbox{\sf Thick}}}
\def\cthi{{\mbox{\sf CThick}}}
\def\late{{\mbox{\sf Late}}}
\def\clate{{\mbox{\sf CLate}}}
\def\plate{{\mbox{\sf PLate}}}
%fraction with round brackets
\newcommand{\ffrac}[2]
  {\left( \frac{#1}{#2} \right)}
\newcommand{\calF}{{\mathcal F}}
\newcommand{\dfn}{\stackrel{\triangle}{=}}
\newcommand{\beqn}[1]{\begin{eqnarray}\label{#1}}
\newcommand{\eeqn}{\end{eqnarray}}
\newcommand{\oo}{\overline}
\newcommand{\uu}{\underline}
\newcommand{\bfcdot}{{\mbox{\boldmath$\cdot$}}}
\newcommand{\Var}{{\rm \,Var\,}}
%qed
\def\squarebox#1{\hbox to #1{\hfill\vbox to #1{\vfill}}}
\renewcommand{\qed}{\hspace*{\fill}
            \vbox{\hrule\hbox{\vrule\squarebox{.667em}\vrule}\hrule}\smallskip}
\newcommand{\half}{\frac{1}{2}\:}
\newcommand{\beaa}{\begin{eqnarray*}}
\newcommand{\eeaa}{\end{eqnarray*}}
\newcommand{\calK}{{\mathcal K}}
\def\dimm{{\overline{{\rm dim}}_{_{\rm M}}}}
\def\ooth{{\hat{\theta}}}
\def\dimp{\dim_{_{\rm P}}}
\def\htaum{{\hat\tau}_m}
\def\htaumk{{\hat\tau}_{m,k}}
\def\htaumkj{{\hat\tau}_{m,k,j}}
\def\loc{{\mbox{\rm loc}}}
\bibliographystyle{alpha}

\title[Tightness for the GFF]
{Tightness of the recentered
maximum of the two--dimensional discrete Gaussian Free Field}

\author[Maury Bramson\,,
 Ofer Zeitouni]
{Maury Bramson$^*$\,\, 
Ofer Zeitouni$^\S$}

\date{September 16, 2010.
\newline\indent
$^*$Research  partially supported by NSF grant number CCF-0729537.
\newline\indent
$^\S$Research  partially supported by NSF grant number
DMS-0804133 and by the Herman P. Taubman chair of Mathematics
at the Weizmann Institute.}

\begin{abstract}
\noindent
We consider the maximum of the discrete two dimensional 
	Gaussian free field (GFF) in a box, and prove that
	its maximum, centered at its mean, is tight, settling a long--standing
	conjecture. The proof 
	combines a recent observation of \cite{BDZ} with elements from
	\cite{bramson2} and comparison theorems for Gaussian fields. 
	An essential part of the argument is the precise evaluation, up to an 
	error of order $1$,
	of the expected value of the maximum of the GFF in a box.
	Related Gaussian fields, such as the GFF on a two--dimensional
	torus, are also discussed.
\end{abstract}

\maketitle

%\ls{2}
\section{Introduction}
\label{sec-introduction}
We consider the discrete Gaussian Free Field (GFF)  in a two-dimensional
box of side $N$, with Dirichlet boundary conditions. That is, let
$V_N=([0,N-1]\cap \mathbb{Z})^2$ and  $V_N^o=((0,N-1)\cap \mathbb{Z})^2$,
and let
$\{w_m\}_{m\geq 0}$ denote a simple random walk started 
in $V_N^o$ and killed at $\tau=\min\{m: w_m\in
\partial V_N\}$ (that is, killed upon hitting the boundary
$\partial V_N=V_N\setminus V_N^o$). 
For $x,y\in V_N$, define
$G_N(x,y)=E^x(\sum_{m=0}^\tau {\bf 1}_{w_m=y})$, where $E^x$ denotes expectation
with respect to the random walk started at $x$. The GFF is
the zero-mean Gaussian field $\{\gff_z^N\}_z$
indexed by $z\in V_N$ with covariance
$G_N$.

Let $\gff_N^*=\max_{z\in V_N} \gff_z^N$.
It was proved in \cite{BDG} that $\gff_N^*/(\log N)\to c$ with $c=2
\sqrt{2/\pi}$;  the proof is closely related to
the proof of the law of large numbers for the
maximal displacement of a branching random walk in $\mathbb{R}$. 

Let $M_N:=\gff_N^*-E\gff_N^*$.
The goal of this paper is to prove the following.
\begin{theorem}
	\label{theo-1}
The sequence of random variables $\{M_N\}_{N\geq 1}$
is tight.
\end{theorem}
The statement in Theorem \ref{theo-1} has been a ``folklore'' conjecture
for some time, and appears in print,
e.g., as  open problem \#4 in \cite{sourav} 
(for an earlier appearance in print of a related conjecture, see
\cite{CLD}).
To the best of our knowledge, prior to the current paper,
the  sharpest result in this direction is due
to \cite{sourav},
 who shows that the variance of $M_N$ is $o(\log N)$, and
to \cite{BDZ}, who show, building on an argument of
\cite{DH91}, that Theorem \ref{theo-1} holds if one replaces $N$ by
an appropriate deterministic sequence $\{N_k\}_{k\geq 1}$.

In the same paper \cite{BDZ}, it is shown that Theorem \ref{theo-1}
holds as soon as one proves that, for an appropriate constant $C$,
$E \gff_{2N}^*\leq  E\gff_N^*+C$ for all $N=2^n$ with $n$ integer. Theorem
\ref{theo-1} thus follows immediately from the following 
theorem, which is our main result.
\begin{theorem}
\label{theo-2}
With notation as above,  
\begin{equation}
\label{eq-main}
E\gff_{2^n}^*= c_1n - c_2  \log n+O(1),
\end{equation}
with $c_1=2\sqrt{2/\pi}\log 2$ and 
%$c_2=3\sqrt{2\log 2/\pi}/2$.
$c_2=(3/4)\sqrt{2/\pi}$.
\end{theorem}
One should note the striking similarity with the behavior of
branching random walks (BRW), see \cite{bramson2} (where branching Brownian 
motions are considered) and \cite{ABR}. The relation
with (an imbedded) BRW is already apparent in \cite{BDG}, but the argument 
there is not sharp enough to allow for a control of the $\log n$ term
in \eqref{eq-main}.

Our approach  
to the proof of Theorem \ref{theo-2} involves two main components. 
The first is a comparison argument (based on the Sudakov-Fernique 
inequality, see Lemma \ref{sud-fer} below), that will allow us to quickly
prove an upper bound in \eqref{eq-main}, and to relate $E\gff_{2^n}^*$ 
to the expectation of the maximum of other Gaussian fields (and,
 in particular,
to a version of the Gaussian Free Field on the torus, denoted 
$\{\tgff_x^N\}_{x\in V_N}$ below, 
as well as to  a modified version of branching random 
walk, denoted $\{\mbrw_x^N\}_{x\in V_N}$ below). 
The second step consists of the analysis of 
the modified branching random walk $\mbrw_x^N$, by properly modifying the
second moment argument in \cite{bramson2} (see also \cite{ABR}). 

Our results also provide
an analog of Theorem \ref{theo-2} for the torus GFF $\{\tgff_x^N\}_{x\in V_N}$,
in Propositions \ref{prop-tgffmbrw} and \ref{prop-lb}, which is of
interest in its own right. Intuitively, the model is the natural
counterpart of the modified branching random walk
$\{\mbrw_x^N\}_{x\in V_N}$, which plays a central role in the proof
of Theorem \ref{theo-2}.  We have not proved
the analog of Theorem \ref{theo-1} for the torus GFF, which requires a
modification of the argument in \cite{BDZ}.

The  paper is structured as follows.
In the next section, we recall a fundamental 
comparison between maxima of 
Gaussian fields; we then
introduce the torus GFF, branching random walk, and 
modified branching random walk, and estimate their covariances. 
Section \ref{sec-ub} is devoted to the proof of the upper bound in
Theorem 
\ref{theo-2}. The rest of the paper deals with the lower bound
in Theorem \ref{theo-2}. Section \ref{sec-lbpre} reduces the proof of
the lower bound to a  lower bound on the maximum
of a truncated version of the modified branching random walk
introduced in Section \ref{sec-prelim}. Section \ref{LBMBRW} reduces the
proof of the latter to a lower bound  on 
the maximum of the modified branching random walk over a subset of $V_N$.
The proof of this bound
is given in Section \ref{sec-6}, using the second moment
method. The proofs of some technical estimates, closely related to estimates
in \cite{bramson2}, are sketched in the appendix.\\

\noindent
{\bf Notation:} throughout, the letter $C$ indicates a positive 
constant, independent
of $N$, whose value may change from line to line. Positive 
constants that are fixed 
once and for all are denoted by the lower case $c$ with a subscript, 
for example $c_5$ or $c_X$.
\section{Preliminaries and  approximations}
\label{sec-prelim}
In this section, we recall a comparison 
tool between the maxima of different Gaussian 
fields and introduce Gaussian
fields that approximate the GFF.

\subsection{The Sudakov--Fernique inequality}
The following inequality allows for the comparison of the expectation
of the maxima of different Gaussian fields. For a proof, see
\cite{sudfer}.
\begin{lemma}[Sudakov--Fernique]
\label{sud-fer}
Let ${\AA}$ denote an arbitrary (finite) set, let 
$\{G_\alpha^i\}_{\alpha\in {\AA}}$, $i=1,2$,
% and $\{g_\alpha\}_{\alpha\in
%{\cal A}}$ 
denote two zero mean Gaussian fields and
set $G^*_i=\max_{\alpha\in {\AA}}G_\alpha^{i}$.
If
\begin{equation}
E(G_\alpha^1-G_\beta^1)^2\geq E(G_\alpha^2-G_\beta^2)^2\,, \quad \mbox{\rm
for all $\alpha,\beta\in {\AA}$}\,,
\end{equation}
then
\begin{equation}
\label{eq-sf}
E G^*_1\geq EG^*_2\,.
\end{equation}
\end{lemma}
In particular, if $\{G_\alpha\}_{\alpha\in {\AA}}$ and $\{g_\alpha\}_{\alpha
\in {\AA}}$ are independent centered
Gaussian fields, then 
one sees that
$$E (\max_{\alpha\in \AA} (G_\alpha+g_\alpha))\geq 
E(\max_{\alpha\in \AA} G_\alpha)\,,$$
a fact that is also easy to check
 without the 
Gaussian assumption.

\subsection{The Torus GFF, Branching Random Walks,
and Modified Branching Random Walks}
We introduce several Gaussian fields with index set $V_N$ that will
play a role in the proof of Theorem
\ref{theo-2}.

\subsubsection{The Torus GFF}
One of the drawbacks of working with the GFF is that its variance is not 
the same at all points of $V_N$. 
The Torus GFF (TGFF) $\{\tgff_z^N\}_{z\in V_N}$
is a Gaussian field whose correlation
structure resembles the GFF, but has the additional property that its variance
is constant across $V_N$. To define it formally, for $x,y\in \Z^2$, write
$x\ \!\! \sim_N \ \!\! y$  if $x-y\in (N\Z)^2$. Similarly, for
$B,B'\subset V_N$, write $B\sim_N B'$ if there exist 
integers $i,j$ so that $B'=B+(iN,jN)$. Let $\tau'$
denote an exponential random variable of parameter $1/N^2$ and, with 
$\{w_m\}_{m\geq 0}$ denoting a simple random walk independent
of $\tau'$, 
define, for $x,y\in V_N,$
$$ \bar G_N(x,y)=E^x(\sum_{m=0}^{\tau'} {\bf 1}_{w_m\ \!\! \sim_N\  \!\!y})\,,$$
where  $E^x$ denotes expectation over both $\tau'$ and
the random walk 
started at $x$.
That is,
$\bar G_N$ is the Green function of a simple random walk on the 
torus of side $N$, killed at the independent exponential time
$\tau'$. The TGFF is the centered Gaussian process
$\{\tgff_z^N\}_{z\in V_N}$ with covariance $\bar G_N$. By construction,
for $x,y\in V_N$, $E((\tgff_x^N)^2)=E((\tgff_y^N)^2)$, and
an easy computation, using known properties of the Green function
of two dimensional simple random walk, see, e.g., \cite{lawler},
reveals that 
\begin{equation}
	\label{eq-vartgff}
	|E((\tgff_x^N)^2)-\frac{2}{\pi} \log N|\leq C \,.
\end{equation}
(Recall that, by our convention on constants, $C$ in \eqref{eq-vartgff}
does not depend on $N$.)
We define $\tgff^*_N=\max_{z\in V_N}\tgff_z^N$.

\subsubsection{Branching Random Walks}
In what follows, we consider $N=2^n$ for some positive
integer $n$. For $k=0,1,\ldots,n$,
let $\BB_k$ denote the collection of subsets of $\Z^2$
consisting of 
squares of side $2^k$ with corners in
$\Z^2$, let $\BB\DD_k$ denote the subset of $\BB_k$ consisting of squares 
of the form $([0,2^k-1]\cap \Z)^2+(i2^k,j2^k)$. Note that the collection
$\BB\DD_k$ partitions $\Z^2$ into disjoint squares. For $x\in V_N$,
let $\BB_k(x)$ denote those elements $B\in \BB_k$ with $x\in B$. Define
similarly $\BB\DD_k(x)$. Note that the set $\BB\DD_k(x)$ contains exactly
one element, whereas $\BB_k(x)$ contains $2^{2k}$ elements.

Let $\{a_{k,B}\}_{k\geq 0, B\in \BB\DD_k}$ denote  an i.i.d. 
family of standard Gaussian random variables.
The BRW $\{\brw_z^N\}_{z\in V_N}$ 
is defined by
$$\brw_z^N=\sum_{k=0}^n \sum_{B\in \BB\DD_k(z)} a_{k,B}\,.$$
We again define $\brw^*_N=\max_{z\in V_N}\brw_z^N$.

\subsubsection{Modified Branching Random Walks}
We continue to consider
$N=2^n$ for some positive
integer $n$ and again employ the notation
$\BB_k$ and $\BB_k(x)$.
Let $\BB_k^N$ denote the collection of subsets of $\Z^2$ consisting of
squares of side $2^k$ with lower left corner in $V_N$.
Let $\{b_{k,B}\}_{k\geq 0, B\in \BB_k^N}$ denote  an i.i.d.
family of centered Gaussian random variables of variance $2^{-2k}$,
and define
$$ b_{k,B}^N=\left\{\begin{array}{ll}
b_{k,B},& B\in \BB_k^N,\\
b_{k,B'},& B\sim_N B'\in \BB_k^N\,.
\end{array}
\right.
$$
The modified 
branching random walk (MBRW) $\{\mbrw_z^N\}_{z\in V_N}$
is defined by
$$\mbrw_z^N=\sum_{k=0}^n \sum_{B\in \BB_k(z)} b_{k,B}^N\,.$$
Note that, by construction, $E((\brw_z^N)^2)
=E((\mbrw_z^N)^2)=n+1$.
We again define $\mbrw^*_N=\max_{z\in V_N}\mbrw_z^N$.

\subsubsection{Geometric distances}
The following are several notions of distances between points in $V_N$. First,
$\|\cdot\|$ denotes the Euclidean norm, while  $\|\cdot\|_\infty$ denotes the
$\ell^\infty$ norm. Thus, for $x,y\in V_N$,
$\|x-y\|$ and $\|x-y|_\infty$ induce metrics with
$$\|x-y\|_\infty\leq \|x-y\|\leq \sqrt{2}\|x-y\|_\infty\,.$$
We also need to consider distances on the torus determined by $V_N$. Those 
are defined by
$$ d^N(x,y)=\min_{z: \ z\  \!\!\sim_N \  \!\!y}\ \|x-z\|\,,\quad
 d^N_\infty(x,y)=\min_{z:\  z\  \!\!\sim_N \  \!\!y} \ \|x-z\|_\infty\,.
$$
\subsection{Covariance comparisons}
We collect in this subsection some basic facts concerning the covariances
of the Gaussian fields introduced earlier. For a centered
Gaussian field
$\{G_z\}$, we write $\RRR_G(x,y)=E(G_xG_y)$ for its covariance function.
Thus, for example, the covariance function of the GFF (on $V_N$) is denoted
by $\RRR_{\gff^N}$. 

The following is an estimate on  
$\RRR_{\tgff^N}$, $\RRR_{\mbrw^N}$ and
$\RRR_{\gff^N}$.

\begin{lemma}
	\label{lem-mbrwcov}
	There exists a constant $C$ so that,
	with $N=2^n$, the following estimates hold:
	for any $x,y\in V_N$,
	\begin{equation}
		\label{eq-comp1}
		|\RRR_{\tgff^N}(x,y)-\frac{2\log 2}{\pi}(n-\log_2 d^N(x,y))|
		\leq C\,
\end{equation}
and
\begin{equation}
\label{eq-comp2}
		|\RRR_{\mbrw^N}(x,y)-(n-\log_2 d^N(x,y))|
		\leq C\,.
	\end{equation}
	Further, for any $x,y\in V_N+(2N,2N)$,
\begin{equation}
		\label{eq-comp1a}
		|\RRR_{\gff^{4N}}(x,y)-\frac{2\log 2}{\pi}(n-(\log_2 \|x-y\|)_+)
		|\leq C\,.
\end{equation}
\end{lemma}
\proof We begin with the estimate \eqref{eq-comp2}
concerning the MBRW. For $x=(x_1,x_2)$ and
$y=(y_1,y_2)$, write, for $i=1,2$,
$t_i(x,y)=\min(|x_i-y_i|,|x_i-y_i-N|,|x_i-y_i+N|)$.
One then has
\begin{eqnarray}
	\RRR_{\mbrw^N}(x,y)&=&
	\sum_{k=\lceil \log_2 (d^N_\infty(x,y)+1)\rceil}^n 2^{-2k}
	\left[2^k-t_1(x,y)\right]\cdot
	\left[2^k-t_2(x,y)\right]\nonumber\\
	&=&
	\sum_{k=\lceil \log_2 (d^N_\infty(x,y)+1)\rceil}^n 
	\left(1-\frac{t_1(x,y)}{2^k}-\frac{t_2(x,y)}{2^k}
	+\frac{t_1(x,y)t_2(x,y)}{4^k}\right)\,.
	\label{eq-280710d}
\end{eqnarray}
Because $a+b-ab\geq 0$ for $0\leq a,b\leq 1$, we get that
\begin{equation}
	\label{eq-240710a}
	\RRR_{\mbrw^N}(x,y)\leq n-
	 \log_2 (d^N_\infty(x,y)+1)+2\leq
	 n-\log_2 (d^N(x,y)+1)+3\,.
 \end{equation}
 On the other hand, using that $a+b-ab\leq a+b$ for $a,b\geq 0$, 
 we get that
\begin{eqnarray}
	\label{eq-240710b}
	\RRR_{\mbrw^N}(x,y)&\geq &n-
	 \log_2 (d^N_\infty(x,y)+1)-
	\sum_{k=\lceil \log_2 (d^N_\infty(x,y)+1)\rceil}^n 
	2^{-(k-1)} d^N_\infty (x,y)\nonumber\\
	&\geq& 
	 n-\log_2 (d^N(x,y)+1)-C\,.
 \end{eqnarray}
 Combining \eqref{eq-240710a} 
 and \eqref{eq-240710b} yields the claimed estimate on
 $\RRR_{\mbrw^N}$.
	 
We next prove the estimate \eqref{eq-comp1a}.
Note that
for $x,y\in V_N+(2N,2N)$,
                \begin{equation}
                        \label{eq-270710d}
                        \RRR_{\gff^{4N}}(x,y)=
                        P^x(\tau_y\leq \tau_{4N})\RRR_{\gff^{4N}}(y,y)\,,
                \end{equation}
                where $\tau_y=\min\{m\geq 0: w_m=y\}$ and
                $\tau_{4N}=\min\{m\geq 0: w_m\not\in V_{4N}\}$.
                Using, e.g., 
\cite[Exercise 1.6.8]{lawler},
                $$      P^x(\tau_y\leq \tau_{4N})=
                \left(1-\frac{(\log \|x-y\|)_+}{n\log 2}\right)+\frac{O(1)}{n}\,.
                $$
                Moreover, for $x\in V_N+(2N,2N)$,
                $$\RRR_{\gff^{4N}}(x,x)=\frac{2\log 2}{\pi}n+O(1)\,,$$
                see, e.g., \cite{sourav}. Combining these estimates
yields \eqref{eq-comp1a}.

 The estimate 
on $\RRR_{\tgff^N}$ 
 in \eqref{eq-comp1}
requires more work but is still straight
 forward. Recall the simple random walk $\{w_m\}$ and, 
 for $y\in V_N$, denote by $[y]_N=\{z\in V_N:\ z\  \!\!\sim_N \  \!\!y\}$
 the collection of points in $\Z^2$ identified with $y$ for the torus.
 Then, by the Markov property and the memoryless property
 of the exponential distribution,
 \begin{eqnarray}
	 \label{eq-240710c}
	 \RRR_{\tgff^N}(x,y)&=&
 E (\tgff^N_y)^2 
 P^x(\{w_m\}\ \mbox{\rm hits
 $[y]_N$ before $\tau'$})\nonumber\\
 &=&
 \frac{2\log 2}{\pi} n
 P^x(\{w_m\}\ \mbox{\rm hits
 $[y]_N$ before $\tau'$})+O(1)\,,
 \end{eqnarray}
 where we recall that $\tau'$ denotes a geometric random variable of mean 
 $N^2$ and we used \eqref{eq-vartgff} in the second equality.

 Let $\eta$ denote the hitting time of the boundary of a (Euclidean)
 ball of radius $N/2$ around $x$, that is
 $$\eta=\min\{m\geq 0: \|w_m-x\|\geq N/2\}\,.$$
 Let $\tau_y$ denote the hitting time of $[y]_N$, that is
 $$\tau_y=\min\{m\geq 0: w_m\in [y]_N\}\,.$$
 Note that the probability in the right side of \eqref{eq-240710c}
 is $P^x(\tau_y<\tau')$. We have
 \begin{eqnarray}
	 \label{eq-260710a}
	 P^x(\tau_y<\tau')&=&P^x(\tau_y<\eta)+
	 P^x(\tau_y<\tau', \tau_y\geq \eta)-
	 P^x(\tau_y<\eta, \tau_y\geq \tau')\nonumber\\
	 &=:&P_1+P_2-P_3\,.
 \end{eqnarray}
 By standard estimates for two dimensional simple random walk,
 see again e.g., \cite[Exercise 1.6.8]{lawler},
 \begin{equation}
\label{p1est}
 |P_1-[n-\log_2|x-y|]_+/n|\leq C/n 
\end{equation}
and, using the memoryless property of the exponential distribution,
$$P_2\leq \max_{z:\  \|z-[y]_N\|\geq N/4} P^z(\tau_y\leq
\tau')\leq C/n\,.$$

To estimate $P_3$, we use the fact (see e.g., \cite[Lemma 10.4]{sourav})
that, for all $m\geq 1$,
\begin{equation}
	\label{eq-250710b}
	P^x(w_m\in [y]_N)\leq
	\left\{\begin{array}{ll}
		\frac{C}{m}e^{-(d^N(x,y))^2/4m}\,,& m\leq N^2,\\
		&\\
		\frac{C}{N^2}\,,& m>N^2\,.
	\end{array}
	\right.
	\end{equation}
	Write $P_m(x,z):=P^x(w_m\in [z]_N)$. Then,
	again using the memoryless property of the exponential distribution
	and the Markov property  of the simple random walk,
	\begin{equation}
		\label{eq-250710c}
		P_3=P^x(\tau'\leq \tau_y<\eta)\leq
	\frac{C}{N^2}\sum_{m=1}^\infty e^{-m/N^2}\sum_{z\in V_N}
	P_m(x,z)P^z(\tau_y<\eta)\,.
\end{equation}

We split the sum in the right side of \eqref{eq-250710c} into three parts, 
according to the range of $m$ in the summation, writing
$P_3=P_{3,1}+P_{3,2}+P_{3,3}$, with the terms in the right side
determined according to 
$m\leq N^2/n$, $m\in (N^2/n, N^2)$ or $m\geq N^2$. We have
$$P_{3,1}=
\frac{C}{N^2}\sum_{m=1}^{N^2/n} e^{-m/N^2}\sum_{z\in V_N}
	P_m(x,z)P^z(\tau_y<\eta)
	\leq 
\frac{C}{N^2}\sum_{m=1}^{N^2/n} e^{-m/N^2}
\leq \frac{C}{n}\,.
$$

Next, consider  $P_{3,3}$: we have, using
\eqref{eq-250710b} in the first inequality and standard 
estimates for simple random walk in the second, see \cite[Exercise 1.6.8]{lawler},
\begin{eqnarray*}
	P_{3,3}&=&
\frac{C}{N^2}\sum_{m=N^2}^{\infty} e^{-m/N^2}\sum_{z\in V_N}
	P_m(x,z)P^z(\tau_y<\eta)\\
	&\leq & \frac{C}{N^4}
\sum_{m=N^2}^{\infty} e^{-m/N^2}\sum_{z\in V_N}
	P^z(\tau_y<\eta)\\
	&\leq & \frac{C}{N^4}
\sum_{m=N^2}^{\infty} e^{-m/N^2}\sum_{z\in V_N}
\left(\frac{n-\log_2 d^N(z,y)}{n}\right)\\
	&\leq & \frac{C}{N^2}
\sum_{z\in V_N}
\left(\frac{n-\log_2 d^N(z,y)}{n}\right)\\
&\leq & \frac{C}{N^2} \sum_{r=1}^n 2^{2r} \left(1-\frac{r}{n}\right)
\\
&=& C\sum_{r=1}^n 2^{2(r-n)}
\left(1-\frac{r}{n}\right)=\frac{C}{n}\sum_{k=1}^n k2^{-2k}\leq \frac{C}{n}
\,.
\end{eqnarray*}

It remains to estimate $P_{3,2}$. Note first that, for $x$ and $y$
fixed and
each integer part
of the value of $d^N(x,z)$, 
%and $d^N(y,z)$, 
there 
are at most $Cr$  many possible points $z\in V_N$ with $d^N(y,z)\in [r,2r]$. 
Also,
due to \eqref{eq-250710b}, we can write
\begin{eqnarray*}
	P_{3,2}&\leq& \frac{C}{n}+
\frac{C}{N^2}\sum_{m=N^2/n}^{N^2} e^{-m/N^2}\sum_{z\in V_N, d^N(x,z)
\leq d_m}
P_m(x,z)P^z(\tau_y<\eta)\\
&=:& \frac{C}{n}+\frac{C}{N^2}\sum_{m=N^2/n}^{N^2}
e^{-m/N^2}
P_{3,2,m}\,,\end{eqnarray*}
	where $d_m=\sqrt{m\log\log m}\wedge \sqrt{2}N$.
	By summing radially (so that the index $k$ runs over the
	possible integer parts of $d^N(x,z)$ and $n-\ell$ runs over
	the possible integer parts of $\log_2 d^N(y,z)$),
	we can estimate $P_{3,2,m}$ (using \eqref{p1est} for the estimate
	$P_m(x,z)\leq Cm^{-1}e^{-k^2/4m}$ and 
\eqref{eq-250710b} for the estimate $P^z(\tau_y<\eta)\leq C \ell/n$)
by
	%\begin{eqnarray*}
	%	&&
	%	\frac{C}{m}
	%	\sum_{k=1}^{d_m} e^{-k^2/4m}
	%	\sum_{\ell=1}^n
%: j\in %[1,N], 
%	%	j+k\geq d^N(x,y)} 
	%	\left(1-\frac{\log_2 j}{n}\right)
	%	\\ &\leq 
	%	&
	$$	\frac{C}{nm}\sum_{k=1}^{d_m} e^{-k^2/4m}
		\sum_{\ell=1}^n \ell 2^{n-\ell} \leq 
		\frac{CN}{nm}\sum_{k=1}^{d_m} e^{-k^2/4m}
		\leq \frac{CN}{n\sqrt{m}}\,.$$
	%\end{eqnarray*}
	Substituting in the expression for $P_{3,2}$, we get
	$$P_{3,2}
	\leq \frac{C}{n}+\frac{CN}{N^2 n}
	\sum_{m=N^2/n}^{N^2} \frac{e^{-m/N^2}}{\sqrt{m}}
	\leq \frac{C}{n}\,.$$
	Combining the estimates on $P_{3,1}, P_{3,2}$ and $P_{3,3}$ and
	substituting in 
		\eqref{eq-250710c} shows that $P_3\leq C/n$.
		Together with \eqref{eq-240710c}, \eqref{eq-260710a} and
		the estimates on $P_1$ and $P_2$, this 
		completes the proof of the claimed
		estimate on $\RRR_{\tgff^N}$ and hence of the lemma.
		\qed
		\section{The upper bound}
		\label{sec-ub}
		Our goal in this section is to provide the upper bound in
		Theorem \ref{theo-2}; this is achieved in Proposition
		\ref{prop-tgffmbrw} below. We begin by
		relating the maxima of the 
		GFF and the TGFF with the MBRW. 
		\begin{lemma}
			\label{lem-compgfftgffmbrw}
			Let $\{g_z\}_{z\in V_N}$ denote a collection of
			i.i.d. standard Gaussian random variables. Then,
			there exists a constant $C_1$ so that
			\begin{equation}
				\label{eq-270710c}
			%	E \gff_N^*
			%	\leq E\gff_N^*+C\,,
			%\quad	
			\max(E \gff_N^*,E\tgff_N^*)
			\leq \sqrt{\frac{2\log 2}{\pi}}
			E(\max_{z\in V_N}(\mbrw_z^N +
				C_1g_z))\,.
			\end{equation}
		\end{lemma}
		\proof We give the argument for the GFF; the argument for
		the TGFF is similar.
		Note first that, by the definitions and an application of
		Lemma \ref{sud-fer},
		\begin{equation}
			\label{eq-270710f}
			E( \gff_N^*)\leq E(\max_{z\in V_N+(2N,2N)} 
			\gff_z^{4N})\,.
		\end{equation}

		On the other hand,
		%for $x,y\in V_N+(2N,2N)$,
		%\begin{equation}
		%	\label{eq-270710d}
		%	\RRR_{\gff^{4N}}(x,y)=
		%	P^x(\tau_y\leq \tau_{4N})\RRR_{\gff^{4N}}(y,y)\,,
		%\end{equation}
		%where again $\tau_y=\min\{m\geq 0: w_m=y\}$ and
		%$\tau_{4N}=\min\{m\geq 0: w_m\not\in V_{4N}\}$.
		%But
		%$$	P^x(\tau_y\leq \tau_{4N})=
		%\left(1-\frac{\log \|x-y\|}{n\log 2}\right)+\frac{O(1)}{n}\,.
		%$$
		%Moreover, for $x\in V_N+(2N,2N)$,
		%$$\RRR_{\gff^{4N}}(x,x)=\frac{2\log 2}{\pi}n+O(1)\,,$$
		%see e.g., \cite{sourav}.
		%Therefore, 
		writing $x_N=x+(2N,2N)$, $y_N=y+(2N,2N)$
		for $x,y\in V_N$ and
		using \eqref{eq-comp1a} of Lemma \ref{lem-mbrwcov},
we have 
		\begin{equation}
			\label{eq-290710c}
			E((\gff^{4N}_{x_N}-\gff^{4N}_{y_N})^2)\leq
		\frac{2\log 2}{\pi}E((\mbrw^{N}_x-\mbrw^{N}_y)^2)+C\,.
	\end{equation}
		Another application of Lemma \ref{sud-fer},
		together with \eqref{eq-270710f}, 
				completes the proof of 
				 \eqref{eq-270710c} for the GFF.
				 \qed

		It follows from \cite{bramson2} and \cite{ABR}
		that
		\begin{equation}
			\label{eq-270710b}
			E \brw_N^*=2\sqrt{\log 2} n -\frac3{4\sqrt{\log 2}}
			\log n+O(1)\,.
		\end{equation}
		%where $c_3=\log 2$ and $c_2=3/2$. 
		(The statement 
		in \cite{bramson2} is given for branching Brownian motions,
		but the argument given there applies to our BRW as well.)
		This fact, together with
		Lemmas \ref{sud-fer} and  
\ref{lem-compgfftgffmbrw},
%and 
%	Lemma 	\ref{sud-fer},
		yields the following upper bound on the GFF and the TGFF.
		\begin{proposition}
			\label{prop-tgffmbrw}
			There exists a constant $C_2$ such that
			\begin{equation}
				\label{eq-270710a}
				\max(E\gff_N^*,E\tgff_N^*)\leq
				2\log 2\sqrt{\frac{2}{\pi}}n-
				\frac34\sqrt{\frac{2}{\pi}}\log n+C_2\,.
			\end{equation}
		\end{proposition}
		\proof
By construction, for $x,y\in V_N$,
$E((\brw_x^N)^2)=E( (\mbrw_x^N)^2)$ and
$\RRR_{\brw^N}(x,y)\leq \RRR_{\mbrw^N}(x,y)+C$.
By Lemma \ref{sud-fer}, this yields the existence of 
a positive integer $\bar C_1$ such that, with $C_1$ and $g_z$ 
as in Lemma
\ref{lem-compgfftgffmbrw},
%$\{g_z\}_{z\in V_N}$ a collection of
%i.i.d. standard Gaussian variables, 
%it holds that
\begin{equation}
	\label{eq-280710a}
	E(\max_{z\in V_N} (\mbrw_z^N+C_1 g_z))\leq
%	\sqrt{\frac{2\log 2}{\pi}}
E(\max_{z\in V_N}(\brw_z^N+\bar C_1g_z))\,.
\end{equation}
Note however that, by construction, 
$$E(\max_{z\in V_N}(\brw_z^N+\bar C_1g_z))\leq 
E(\brw_{2^{ {\bar C_1}^2} N}^*)\,.$$
Combining this with \eqref{eq-270710b}, 
\eqref{eq-280710a} and
Lemma 
\ref{lem-compgfftgffmbrw} completes the proof of 
\eqref{eq-270710a} and hence of 
Proposition \ref{prop-tgffmbrw}.
\qed

\section{The lower bound: preliminaries}
\label{sec-lbpre}
In this section, we bound from below
the expected maxima of the GFF and the TGFF by an appropriate 
truncation of the MBRW. An analysis of the latter, provided in 
Sections \ref{LBMBRW} and \ref{sec-6}, will then complete the proof of
Theorem \ref{theo-2}.

We begin by introducing the truncation of the 
MBRW alluded to above. Recall that
$$\mbrw_z^N=\sum_{k=0}^n \sum_{B\in \BB_k(z)} b_{k,B}^N\,.$$
For a non-negative integer $k_0\leq n$, define
$$\mbrw_z^{N,k_0}=\sum_{k=k_0}^n \sum_{B\in \BB_k(z)} b_{k,B}^N\,,$$
and write
$\mbrw_{N,k_0}^*=\max_{z\in V_N} \mbrw_z^{N,k_0}$.
Clearly, $\mbrw^{N,0}=\mbrw^N$.
Define, for $x,y\in V_N$,
$\rho_{N,k_0}(x,y)=
E( (\mbrw_x^{N,k_0}-\mbrw_y^{N,k_0})^2)$.
The following are basic properties of
$\rho_{N,k_0}$.
%
%Note that for $x,y\in V_N$,
%$\chi_{N,k_0}(x,y)=
%E( (\mbrw_x^{N,k_0}-\mbrw_y^{N,k_0})^2)$ decreases in $k_0$.
\begin{lemma}
	\label{lem-k0}
	The function $\rho_{N,k_0}$ has the following properties.\\
	\begin{eqnarray}
	&&\rho_{N,k_0}(x,y)\ \mbox{\rm decreases in $k_0$}.\label{k0prop1}\\
	&& 
	\limsup_{k_0\to\infty}
	\limsup_{N\to\infty} 
	\sup_{x,y\in V_N: d_N(x,y)\leq 2^{\sqrt{k_0}}}
	\rho_{N,k_0}(x,y)=0\,. \label{k0prop2}\\
	&&\mbox{\rm There is a function $g:\Z_+\to \R_+$ so that
	$g(k_0)\to_{k_0\to \infty} \infty$ }
	 \nonumber\\
	&& \mbox{\rm and, for
	$x,y\in V_N$ with $d^N(x,y)\geq 2^{\sqrt{k_0}}$,}
	\label{k0prop3}
	\\
	&& \rho_{N,k_0}(x,y)\leq \rho_{N,0}(x,y)-g(k_0)\,, \quad n>k_0.
	\nonumber
\end{eqnarray}
\end{lemma}
\proof 
As in 
	\eqref{eq-280710d} and employing the same notation, 
	we have,
	 for $x\neq y$,
\begin{eqnarray}
	\label{eq-280710h}
	&&\rho_{N,k_0}(x,y)\\
	&=&
	\sum_{k=\lceil \log_2 (d^N_\infty(x,y)+1)\rceil\vee k_0}^n 
	2\left(\frac{t_1(x,y)}{2^k}+\frac{t_2(x,y)}{2^k}
	-\frac{t_1(x,y)t_2(x,y)}{4^k}\right)\nonumber\\
	&&\quad+2(\lceil \log_2 (d_\infty^N(x,y)+1)
	\rceil-k_0)_+
	\,.
	\nonumber
\end{eqnarray}
All properties follow at once from this representation. Indeed, \eqref{k0prop1}
and \eqref{k0prop2} are immediate whereas, to see \eqref{k0prop3}, note that,
for $ \log_2 d_\infty^N(x,y)\geq  \sqrt{k_0}-1$, 
$$\rho_{N,0}(x,y)-\rho_{N,k_0}(x,y)\geq  \sqrt{k_0}-1\,.
$$
\qed

An immediate corollary of Lemmas \ref{sud-fer}, \ref{lem-mbrwcov}
and \ref{lem-k0}
is the following domination by the TGFF
of a truncated MBRW.
\begin{corollary}
	\label{cor-tgffmbrw}
	There exists a constant $k_0$ such that, for all $N=2^n$ large
	and all $x,y\in V_N$,
	\begin{equation}
		\label{eq-290710a}
		\frac{2\log 2}{\pi}
		\rho_{N,k_0}(x,y)\leq E( (\tgff_x^N-\tgff_y^N)^2)\,.
	\end{equation}
	In particular,
	\begin{equation}
		\label{eq-290710aa}
		E\tgff_N^*\geq \sqrt{\frac{2\log 2}{\pi}}E\mbrw_{N,k_0}^*\,.
	\end{equation}
	\end{corollary}

We also need a comparison between the maxima of the GFF and of the MBRW.
	Note that
	\begin{equation}
		\label{eq-290710b}
		\gff_N^*\geq \max_{z\in V_{N/4}+(N/2,N/2)} \gff_z^N\,.
	\end{equation}
	On the other hand, for $x,y\in V_{N/4}$, we have, 
	with the same proof as that of \eqref{eq-290710c},
\begin{equation}
			\label{eq-290710d}
			E((\gff^{N}_{x+(N/2,N/2)}-\gff^{N}_{y+(N/2,N/2)})^2)\geq
		\frac{2\log 2}{\pi}E((\mbrw^{N/4}_x-\mbrw^{N/4}_y)^2)-C\,.
	\end{equation}
Using again
Lemmas \ref{sud-fer}
and \ref{lem-k0}, we get
the following domination by the GFF
of a truncated MBRW.
\begin{corollary}
	\label{cor-gffmbrw}
	There exists a constant $k_0$ such that, for all $N=2^n$ large,
	and all $x,y\in V_N$,
	\begin{equation}
		\label{eq-290710ab}
		E\gff_N^*\geq \sqrt{\frac{2\log 2}{\pi}}E\mbrw_{N/4,k_0}^*\,.
	\end{equation}
\end{corollary}

\section{A lower bound for the truncated MBRW}
\label{LBMBRW}
In this section, we present the proof of the lower bound in Theorem
\ref{theo-2} and an analogous bound for the TGFF. 
That is, we prove the following.
\begin{proposition}
	\label{prop-lb}
The following holds:
\begin{equation}
\label{eq-mainilb}
E\gff_{2^n}^*\geq  c_1n - c_2  \log n+O(1)\,,\quad
E\tgff_{2^n}^*\geq  c_1n - c_2  \log n+O(1)\,,
\end{equation}
with $c_1=2\sqrt{2/\pi}\log 2$ and $c_2=(3/4)\sqrt{2/\pi}$.
\end{proposition}
In view of Corollaries \ref{cor-tgffmbrw} and \ref{cor-gffmbrw},
it is immediate that 
Proposition \ref{prop-lb} follows from the following 
proposition, whose proof will take the rest of this section and the next one.

\begin{proposition}
	\label{proplbmbrw}
	There exists a function $f:\Z_+\to \R_+$ such that,
	for all $N\geq 2^{2 k_0}$,
	\begin{equation}
		\label{eq-290710e}
		E\mbrw_{N,k_0}^*\geq (2 \sqrt{\log 2})n-(3/(4\sqrt{\log 2}))
\log n-f(k_0)\,.
	\end{equation}
\end{proposition}
When proving Proposition \ref{proplbmbrw}, it will be more convenient to 
restrict the maximum in the definition of $S^*_{N,k_0}$ to a subset of $V_N$.
Toward this end, set $V_N'=V_{N/2}+(N/4,N/4)\subset V_N$ and define
$$\tilde S^*_{N,k_0}=\max_{z\in V_N'} S_z^{N,k_0}\,,\quad
\tilde S^*_N=\tilde S^*_{N,0}\,.$$
The  main ingredient in the proof of Proposition \ref{proplbmbrw}
is a lower bound on the upper tail of the distribution
of $\tilde \mbrw_{N,0}^*$, as given in Proposition \ref{proplbmbrw1}
below. 
\begin{proposition}
\label{proplbmbrw1}
Let $A_n=(2\sqrt{\log 2})n-(3/(4\sqrt{\log 2}))\log n$.
There exists a constant $\delta_0\in (0,1)$
%function $\delta_0:\R_+\to (0,1)$
such that, for all $N$,
	\begin{equation}
		\label{eq-300710a}
		P(\tilde \mbrw_{N}^*\geq A_n)\geq \delta_0\,.
	\end{equation}
\end{proposition}
The proof of Proposition \ref{proplbmbrw1} is technically involved
and is deferred to Section \ref{sec-6}.
In the rest of this section, we
show how Proposition \ref{proplbmbrw} follows from
Proposition \ref{proplbmbrw1}.\\

\noindent
{\it Proof of Proposition \ref{proplbmbrw} (assuming
Proposition \ref{proplbmbrw1})}
Our plan is to show that the left tail of $\tilde S_{N}^{*}$ is decreasing
exponentially fast; together with the bound \eqref{eq-300710a}, 
this will imply \eqref{eq-290710e}
with $k_0=0$.  At the end of the proof, we show how the bound for $k_0>0$
follows from the case $k_0=0$.  In order to show the exponential decay, we
compare $\tilde{S}_{N}^{*}$, after appropriate truncation, to four
independent copies of the maximum over smaller boxes, and then iterate.

%We first consider $k_0=0$, and, at the end of the proof,
%we show how the bound for
%$k_0>0$ follows from the case $k_0=0$.

For $i=1,2,3,4$, introduce the four sets 
% $W_{N,i}=[N/8,3N/8)^2+z_i$ where $z_1=(0,0)$, $z_2=(N/2,0)$,
 $W_{N,i}=[0,N/32)^2+z_i$ where $z_1=(N/4,N/4)$, $z_2=(23N/32,N/4)$,
$z_3=(N/4,23N/32)$ and $z_4=(23N/32,23N/32)$.
(We have used here that $3/4-1/32=
23/32$.) 
Note that $\cup_i W_{N,i} \subset V_N$, and that these sets are $N/4$-separated,
that is, for $i\neq j$,
$$\min_{z\in W_{N,i}, z'\in W_{N,j}} d_\infty^N(x,y)> N/4\,.$$ 
Recall the definition of
$\mbrw_z^{N}$ and define, for $n>6$,
$$\bar\mbrw_z^{N}=\sum_{k=0}^{n-6} \sum_{B\in \BB_k(z)} b_{k,B}^N\,;$$
note that
%, for $n>k_0+6$,
$$\mbrw_z^{N}-\bar\mbrw_z^{N}=
\sum_{j=0}^5
 \sum_{B\in \BB_{n-j}(z)} b_{n-j,B}^N
\,.$$
%We have that 

Our first task is to bound the probability
that 
$\max_{z\in V_N}(\mbrw_z^{N}-\bar\mbrw_z^{N})$ is large.
This will be achieved by applying Fernique's criterion 
in conjunction with Borell's inequality. We introduce
some notation. Let $m(\cdot)=m_N(\cdot)$ denote the uniform probability
measure on $V_N$ (i.e., the counting measure normalized by $|V_N|$)
and let  $g:(0,1]\to \R_+$ be the function defined by
$$ g(t)=\left( \log (1/t)\right)^{1/2}\,.$$
Set $\G_z^{N}
=\mbrw_z^{N}-\bar\mbrw_z^{N}$
%In the notation of Lemma 
%\ref{lem-k0}, we have that for $n> k_0+6$,
%$$ E(\G_z-G_{z'})^2=\rho_{N,2^{-5}N}(z,z')\,.$$
and
$$B(z,\epsilon)=\{z'\in V_N: E( (\G_z^N-\G_{z'}^N)^2)
\leq \epsilon^2\}\,.$$
Then,
Fernique's criterion,
see \cite[Theorem 4.1]{adler}, 
implies that, for some universal
 constant $K\in (1,\infty)$, 
\begin{equation}
\label{eq-fer} E(\max_{z\in V_N}\G_z^N)\leq
K\sup_{z\in V_N} \int_0^\infty g(m(B(z,\epsilon))) d\epsilon\,.
\end{equation}

For $n\geq 6$, we have,
in the notation of Lemma 
\ref{lem-k0}, 
$$ E((\G_z^N-\G_{z'}^N)^2)=\rho_{N,n-5}(z,z')\,.$$
Therefore, employing
\eqref{eq-280710h}, there exists a constant
$C$ such that, for $\epsilon\geq 0$, 
$$ \{z'\in V_N: d^N_\infty(z,z')\leq \epsilon^2 N/C\}\subset B(z,\epsilon)\,.$$
In particular, for $z\in V_N$ and $\epsilon>0$,
$$ m(B(z,\epsilon))\geq 
( (\epsilon^4 /C^2)\vee (1/N^2) )\wedge 1\,.$$ 
Consequently,
$$\int_0^\infty g(m(B(z,\epsilon))) d\epsilon
\leq \int_0^{\sqrt{C/N}} \sqrt{\log(N^2)} d\epsilon+
\int_{\sqrt{C/N}}^{\sqrt{C}} \sqrt{\log(C^2/\epsilon^4)} d\epsilon<C_4\,,$$
for some constant $C_4$. Applying Fernique's criterion
\eqref{eq-fer}, we deduce that
$$E (\max_{z\in V_N}
(\mbrw_z^{N}-\bar\mbrw_z^{N}))\leq C_4K\,.$$

The expectation
$E((\mbrw_z^{N}-\bar\mbrw_z^{N})^2)$ is  bounded in
$N$.
%and $k_0$. 
Therefore, using Borell's inequality, see, e.g.,
\cite[Theorem 2.1]{adler}, it follows that,
for some constant $C_5$ and all $\beta>0$,
\begin{equation}
\label{eq-300710b}
P( \max_{z\in V_N}
(\mbrw_z^{N}-\bar\mbrw_z^{N})\geq 
C_4K+\beta)\leq 2e^{-C_5 \beta^2}\,.
\end{equation}

We also note the following bound, which is obtained
similarly: there exist constants $C_5,C_6$ such that,
for all $\beta>0$,
\begin{equation}
\label{eq-300710bb}
P( \max_{z\in V_{N/16}'}
(\bar\mbrw_z^{N}-\mbrw_z^{N/16})\geq 
C_6+\beta)\leq 2e^{-C_7 \beta^2}\,.
\end{equation}

The advantage of working with $\bar \mbrw^{N}$ instead of
$\mbrw^{N}$ is that the fields
$\{\bar\mbrw_z^{N}\}_{z\in W_{N,i}}$ are independent for $i=1,\ldots,4$.
For every $\alpha,\beta> 0$,
we have the bound
\begin{eqnarray}
\label{eq-300710e}
&&P(\tilde \mbrw_{N}^*\geq A_n-\alpha)\\
&\geq &
P(\max_{z\in V_N'}
\bar\mbrw^{N}_z\geq A_n+C_4-\alpha+\beta)-
P( \max_{z\in V_N'}
(\mbrw_z^{N}-\bar\mbrw_z^{N})\geq 
C_4+\beta)\nonumber\\
&\geq&
P(\max_{z\in V_N'}\bar\mbrw_{N}\geq A_n+C_4-\alpha+\beta)-
2e^{-C_5 \beta^2}\,,
\nonumber
\end{eqnarray}
where \eqref{eq-300710b} was used in the last inequality.
On the other hand, for any $\gamma,\gamma'>0$,
\begin{eqnarray*}
&&P(\max_{z\in V_N'}\bar\mbrw_z^{N}\geq A_n-\gamma)\geq
P(\max_{i=1}^4\max_{z\in W_{N,i}}
\bar\mbrw_z^{N}
\geq A_n-\gamma)\\
&=&1-(
P(\max_{z\in W_{N,1}}
\bar\mbrw_z^{N}
<A_n-\gamma))^4\\
%&\geq &
%1-( P(\max_{z\in W_{N,1}} \mbrw_z^{N,k_0}
%<A_n-\gamma+C_4+\gamma')+2e^{-C_5(\gamma')^2})^4\\
&\geq&
1-\left(
P(\max_{z\in V_{N/16}'}
\mbrw_z^{N/16}
<A_n-\gamma+C_6+\gamma')+2e^{-C_7(\gamma')^2}\right)^4
\,,
\end{eqnarray*}
where \eqref{eq-300710bb} was used in the inequality.
Combining this estimate with \eqref{eq-300710e}, we get that,
for any $\alpha,\beta,\gamma'>0$,
\begin{eqnarray}
\label{eq-300710f}
&&P(\tilde \mbrw_{N}^*\geq A_n-\alpha)
\\
&\geq& 1-2e^{-C_5 \beta^2}\nonumber\\
&&
\!\!\!\!\!\!
\!\!\!\!\!\!
\!\!\!\!\!\!
-\left(
P(\max_{z\in V_{N/16}'}
\mbrw_z^{N/16}
<A_n+C_4+C_6+\beta+\gamma'-\alpha)+2e^{-C_7(\gamma')^2}\right)^4\,.
\nonumber
\end{eqnarray}

We now iterate the last estimate. Let $\eta_0=1-\delta_0<1$
and, for $j\geq 1$, choose a constant $C_8=C_8(\delta_0)>0$
so that, for 
$\beta_j=\gamma_j'=C_8 \sqrt{\log (1/\eta_j)}$,
%for appropriate $C_8$.)
%$\beta_j,\gamma'_j>0$  so 
%that 
$$ \eta_{j+1}=2e^{-C_5 \beta_j^2}+(\eta_j+2e^{-C_7 (\gamma_j')^2})^4$$
satisfies  $\eta_{j+1}<\eta_j(1-\delta_0)$.
(It is not hard to verify that such a choice is possible.)
%(When $\eta_j$ is small enough, say smaller than $1/4$, choose
%$\beta_j=\gamma_j'=C_8 \sqrt{\log (1/\eta_j)}$ for appropriate $C_8$.)
With this choice of $\beta_j$ and $\gamma_j'$,
 set $\alpha_0=0$ and $\alpha_{j+1}=\alpha_j+C_4+C_6+\beta_j+
\gamma_j'$, noting that
$\alpha_j\leq C_9 \sqrt{\log(1/\eta_j)}$ for some $C_9=C_9(\delta_0)$. 
Substituting in
\eqref{eq-300710f} and using Proposition 
\ref{proplbmbrw} to start the recursion,
we get that
\begin{equation}
\label{eq-020910a}
%P(\mbrw_{N}^*\geq A_n-\alpha_{j+1})\geq
P(\tilde \mbrw_{N}^*\geq A_n-\alpha_{j+1})
\geq 1-\eta_{j+1}\,.
\end{equation}
Therefore, 
%using the last estimate with $y=0$,
\begin{eqnarray*}
%E\mbrw_N^*&\geq& 
E\tilde\mbrw_{N}^*&\geq& 
A_n-\int_0^\infty
P(\tilde\mbrw_{N}^*\leq A_n-\theta)d\theta\\
%&\geq& A_n-\sum_{j=0}^\infty ({\alpha_j}-\alpha_{j+1})
&\geq& A_n-\sum_{j=0}^\infty \alpha_j
P(\tilde\mbrw_{N}^*\leq A_n-\alpha_j)\\
&\geq &A_n-
C_9 \sum_{j=0}^\infty \eta_j\sqrt{\log(1/\eta_j)}\,.
\end{eqnarray*}
Since $\eta_j\leq (1-\delta_0)^j$, it follows that there exists a constant
$C_{10}>0$ so that  
\begin{equation}
	\label{eq-020910g}
E\mbrw_{N}^*
\geq E\tilde \mbrw_{N}^*
\geq A_n-C_{10}\,.
\end{equation}
This completes the proof of Proposition \ref{proplbmbrw} in the
case $k_0=0$.

To consider the case $k_0>0$, define 
$$\hat S^*_{N,k_0}=\max_{z\in V_N'\cap 2^{k_0}\Z^2} \mbrw_z^{N,k_0}\,.
$$
Then, $\hat S^*_{N,k_0}\leq \tilde S^*_{N,k_0}$. On the other hand,
$\hat S^*_{N,k_0}$ has, by construction, the same distribution 
as $\tilde S^*_{2^{-k_0}N,0}=\tilde S^*_{2^{-k_0}N}$. Therefore,
for any $y\in \R$,
%\begin{eqnarray*}
$$ P(\tilde S^*_{N,k_0}\geq y)\geq
P(\hat S^*_{N,k_0}\geq y)
\geq
P(\tilde S^*_{2^{-k_0}N}\geq y)\,.
$$
We conclude that
$$
ES^*_{N,k_0}\geq E\tilde S^*_{N,k_0}\geq
E\tilde S^*_{2^{-k_0}N}\,.$$
Application of \eqref{eq-020910g} completes the proof 
of Proposition \ref{proplbmbrw}.
\qed

\section{Proof of Proposition \ref{proplbmbrw1}}
\label{sec-6}
The proof is based on the second moment method and is very similar
to the argument in \cite[Section 6]{bramson2}. 
%We first reduce to the case
%$k_0=0$. Indeed, define
%$$\hat S^*_{N,k_0}:=\max_{z\in V_N'\cap 2^{k_0}\Z^2} \mbrw_z^N\,.
%$$
%Then, $\hat S^*_{N,k_0}\leq \tilde S^*_{N,k_0}$.
%
%
%For ease of notation, we 
%work with $k_0=0$ throughout, but the argument applies equally well (with 
%all constants depending on $k_0$) in the general case as well, if one 
%uses it for  the variable 
%$\max_{z\in V_N'\cap 2^{k_0}\Z^2} \mbrw_z^N\leq \tilde S^*_{N,k_0}$.
We begin by introducing some
notation. Recall that, for $z\in V_N$, 
$$ \mbrw_z^{N}=\sum_{k=0}^n \sum_{B\in 
\BB_k(z)} b_{k,B}^N\,.$$
We introduce  a time parameter, setting 
\begin{equation}
\label{eq-mbrwtime}
 \mbrw_z^N(j)=\sum_{k=n-j}^n \sum_{B\in 
\BB_k(z)} b_{k,B}^N\,
\end{equation}
for $j=0,\ldots,n$.
Fix a (large) constant $c_5$ and introduce the 
function $L_n(j)$, $j=0,1,\ldots,n$, with $L_n(0)=L_n(n)=0$ and
$$ L_n(j)=\left\{\begin{array}{ll}
c_5 \log j,& j=1,\ldots,\lfloor n/2\rfloor\\
c_5 \log (n-j),& j=\lfloor n/2\rfloor+1,\ldots,n-1\,.
\end{array}
\right. $$
	
We next introduce events involving the path $\mbrw_z^N(\cdot)$.
Recall that $A_n=(2\sqrt{\log 2})n-  (3/(4\sqrt{\log 2}))\log n$.
For $z\in V_n$, define the event
$$ \CC_z=\{\mbrw_z^N(j)\leq \frac{j}{n} (A_n+1)-L_n(j)+1, j=0,1,\ldots,n,
\mbrw_z^N\in [A_n,A_n+1]\}\,.$$	
Define
$$ h=\sum_{z\in V_N'} {\bf 1}_{\CC_z}\,.$$	
We have the following proposition.
\begin{proposition}
\label{prop-m1}
There exists a constant $C_{11}>0$ such that
\begin{equation}
\label{eq-260810a}
Eh\geq C_{11}^{-1}
\end{equation}
and
\begin{equation}
\label{eq-260810b}
Eh^2\leq C_{11}
\end{equation}
\end{proposition}
Proposition \ref{proplbmbrw1} 
%(for $k_0=0$)
follows at once from Proposition
\ref{prop-m1}, the Cauchy-Schwartz inequality $P(h\geq 1)\geq
(Eh)^2/E(h^2)$,
and the inequality
$$ P(\tilde S^*_{N}\geq A_n)\geq P(h\geq 1)\,.$$
%As mentioned before, the same argument applies to any fixed $k_0$.
The rest of this section is devoted to the proof of Proposition
\ref{prop-m1}. In the proof, certain crucial estimates (Lemmas
%equation \eqref{eq-260810c} and Lemma 
\ref{lem-mauryinsist} and \ref{lem-11}) are discrete analogues
of corresponding results in \cite{bramson2}. We provide
in the appendix some detail on the proof of these lemmas.\\

\noindent
{\bf Proof of Proposition \ref{prop-m1}}
We begin with the following lemma;
%, whose proof
%parallels that of
%by noting that, with the same  proof as in 
%\cite[Lemma 9]{bramson2}; 
the appendix supplies details on the proof.
%see the appendix for 
\begin{lemma}
	\label{lem-mauryinsist}
For some $C_{12},C_{13}>0$,
\begin{equation}
\label{eq-260810c}
C_{12}N^{-2}\geq P(\CC_z)\geq C_{13} N^{-2}\,.
\end{equation}		
\end{lemma}
It follows from Lemma \ref{lem-mauryinsist} that
%Therefore
$$ Eh\geq  C_{12} |V_N'|/|V_N|= C_{12}/4\,,$$
proving \eqref{eq-260810a}.
%(A crucial aspect here is that $P(A_z)$ remains the same if the process
%$\brw_z^N$ would be substituted for $\mbrw_z^N$.)
		
To compute the second moment in \eqref{eq-260810b}, we first set
%introduce some notation.
$$r(z,z')=n-\lceil \log_2 (d_N^\infty(z,z')+1)\rceil\,,$$
for $z,z'\in V_N'$.
%% r is the amount of time they go together
A crucial observation is that, by construction, 
\begin{eqnarray}
	\label{eq-270810a}
&& \mbox{\rm the process 
$\{\mbrw_{z'}^N(\ell+r(z,z'))-\mbrw_{z'}^N(r(z,z'))\}_{\ell\geq 0}$}\nonumber\\
&&
\mbox{\rm is independent of  the sigma algebra generated by}\\
&&\mbox{\rm the processes $\{\mbrw_z^N(j)\}_{j\geq 0}$ and 
$\{\mbrw_{z'}^N(j)\}_{j\leq r(z,z')}$.}\nonumber
\end{eqnarray}
(Note that the boxes involved in the construction of the first process are
disjoint from those of the other two processes.)
We employ the decomposition
%\begin{eqnarray*}
%	Eh^2& =&\sum_{z,z'\in V_N'} P(\CC_z, \CC_{z'})\\
%&=&
%\sum_{z,z': r(z,z')<n/2} P(\CC_z, \CC_{z'})
%+
%\sum_{z,z': r(z,z')>n/2} P(\CC_z, \CC_{z'})=:Q_1+Q_2\,.
%\end{eqnarray*}
\begin{eqnarray*}
	Eh^2 &=&\sum_{z,z'\in V_N'} P(\CC_z, \CC_{z'})\\
	&=&
\sum_{z,z': r(z,z')<n/2} P(\CC_z, \CC_{z'})
%+
%\sum_{z,z': r(z,z')\in [n/6,n/2]} P(\CC_z, \CC_{z'})\\
%&& +
+\sum_{z,z': r(z,z')\geq n/2} P(\CC_z, \CC_{z'})
=: Q_1+Q_2\,.\end{eqnarray*}

We begin by considering $Q_2$. For this,
we introduce the event
$$ \tilde \CC_{z,z'}=\{
\mbrw_{z'}^N(r(z,z'))\leq \frac{r(z,z')}{n} (A_n+1)-L_n(r(z,z'))+1, 
\mbrw_{z'}^N\in [A_n,A_n+1]\}\,,$$	
noting that $\CC_{z'}\subset \CC_{z,z'}$. 
	It follows from \eqref{eq-270810a} that
	\begin{eqnarray}
		\label{eq-270810b}
		P(\CC_z,\CC_{z'})&\leq &
		P(\CC_z,\tilde \CC_{z,z'})\nonumber\\
	&\leq&
		P(\CC_z)P\left(G_{z,z'}
		\geq \left(1-\frac{r(z,z')}{n}\right)(A_n+1)+L_n(r(z,z'))\right)\,,
	\end{eqnarray}
	where $G_{z,z'}$ is a centered Gaussian random variable of variance
	$$n-r(z,z')=\lceil \log_2 (d_N^\infty(z,z')+1)\rceil=: u(z,z').$$
	Therefore,
	using \eqref{eq-260810c},
	\begin{eqnarray}
		\label{eq-270810d}
	\nonumber	P(\CC_z,\CC_{z'})&\leq& 
		C_{14}N^{-2} 
		\exp\left(-((A_n/n)u(z,z')+ L_n(r(z,z')))^2/2 u(z,z')\right)\\
		&\leq &C_{15} 2^{-2n - 2\log_2
		%d_N^{\infty}(z,z')} e^{c_6\log n}
		d_N^{\infty}(z,z')} e^{3 \log n}
		e^{-(A_n/n)L_n(r(z,z'))}\,.
	\end{eqnarray}
	%where the constant $c_6$
	%=3/4\sqrt{\log 2}$ 
	%does not depend
	%on $c_5$.
	Since the number of points $z'$ with $d_N^\infty(z,z')\in [2^k,2^{k+1}]$
	is bounded by a constant multiple of $2^{2k}$, we conclude from
	\eqref{eq-270810d} that
	$$Q_2 \leq
	C_{16}\sum_{k=n/2}^n (n-k)^{-c_5A_n/n}n^{3}<C_{17}\,,$$
	if $c_5$ is chosen large enough.
	%Similarly,
	%$$Q_2 \leq
	%C_{16}\sum_{k=n/6}^{n/2} k^{-c_5(A_n/n)}n^{c_6}<C_{18}\,.$$

	It thus remains to handle  $Q_1$. 
Introduce the events	
$$\DD_{z,z'}^{(1)}=\{
\mbrw_{z'}^N(r(z,z'))\leq \frac{r(z,z')}{n} (A_n+1)-L_n(r(z,z'))+1 \}$$
and, for $w\in \R$, 
\begin{eqnarray*}
	 \DD_{z,z',w}^{(2)}&=&\{
	\mbrw_{z'}^N(j)-\mbrw_{z'}^N(r(z,z'))\leq 
	\frac{j}{n} (A_n+1)-w+1, j=r(z,z'),\ldots,n,\\
	&&\mbrw_{z'}^N(n)-\mbrw_{z'}^N(r(z,z'))\in [A_n-w,A_n+1-w]\}\,.
\end{eqnarray*}
	It follows again from
	\eqref{eq-270810a} that
	\begin{eqnarray}
		\label{eq-270810bnew}
		P(\CC_z,\CC_{z'})&\leq &
		P(\CC_z,\DD_{z,z'}^{(1)},\DD_{z,z',\mbrw_{z'}^N(r(z,z'))}^{(2)})\nonumber\\
	&\leq&
	P(\CC_z)\max_{w\leq r(z,z') (A_n+1)/n-L_n(r(z,z'))+1}
	P(\DD_{z,z',w}^{(2)})\,.
\end{eqnarray}

To analyze 
	$P(\DD_{z,z',w}^{(2)})$,
we employ the following lemma.
(Details of the proof are given 
in the appendix.) 
%, whose proof mirrors that of
%\cite[Lemma 11]{bramson2}.
\begin{lemma}
	\label{lem-11}
	With notation as above, there exist  constants $C_{19}$ and
	$C_{20}$  
	so that, if $r(x,x')\leq n/2$ and 
	$w\leq 
	r(z,z') (A_n+1)/n-L_n(r(z,z'))+1$, then
	\begin{equation}
		\label{eq-270810f}
		P(\DD_{z,z',w}^{(2)})\leq
		C_{20} (L(r(z,z'))+1)
		\cdot  2^{-2\log_2d_N^\infty(z,z')} \cdot 
		e^{- C_{19}L(r(z,z'))}\,.
	\end{equation}
\end{lemma}
%(We provide a sketch of the proof of Lemma \ref{lem-11}
%in the appendix.) 
%Using Lemma \ref{lem-11} and 
Substituting  (\ref{eq-270810f}) of the lemma 
into 
(\ref{eq-270810bnew}),
 we conclude that 
$$ Q_1\leq C_{20} \sum_{k=0}^{n/2} (L(k)+1)e^{-C_{19} L(k)}\leq C_{21}\,.$$
This completes the proof of Proposition \ref{prop-m1}.
\qed
	
\section{Appendix}	

\begin{appendix}
\setcounter{equation}{0}
\renewcommand{\theequation}{A.\arabic{equation}}
In this appendix, we provide more detail on the bounds in (\ref{eq-260810c})
of Lemma \ref{lem-mauryinsist}
%Proposition 6.2 
and \eqref{eq-270810f} of Lemma \ref{lem-11}, 
%of Proposition 6.3, 
which in both cases are very similar to
material in \cite{bramson2}.

%Bounds in Proposition 6.2.  
\subsubsection*{Bounds in Lemma \ref{lem-mauryinsist}}
%We begin with the bounds in (\ref{eq-260810c}).
Using the Brownian bridge $\frak{z}(s)$, $s\in
[0,n]$, that is standard Brownian motion $\frak{x}(s)$, $s\ge 0$,
conditioned on $\frak{x}(n)=0$, one has
\begin{equation}
\label{eqA1}
P(\frak{z}(s)<2, s=0,\ldots,n) \ge P(\mathcal{C}_z) / K(n)
\ge P(\frak{z}(s) < 1-L(s), s=0,\ldots,n),
\end{equation}
where $K(n) = P(\frak{x}(n) \in [A_n, A_n+1]$.
%
%**You should likely change the symbols here.  I tried to choose something
%similar to my thesis, but do not know how these will look.  The numbers of
%succeeding constants will of course need to be changed.**
%
One can check that $(N^2 / n)K(n) \in [C_{22},C_{23}]$, for 
$0 < C_{22}< C_{23} <
\infty$.  So, in order to demonstrate Lemma \ref{lem-mauryinsist},
it suffices to show the
bounds on each side of (\ref{eqA1}) are of order $1/n$.

By \cite[Proposition 2$^\prime$ on page 555]{bramson2},
% of [Br78],
%
\begin{equation}
\label{eqA2} P(\frak{z}(s) < 1-L(s), s\in [0,n]) \ge C_{24} /n,
\end{equation}
which gives the desired lower bound.  One obtains the analogous upper bound
$C_{25} /n$ for the left side
of (\ref{eqA1}) by applying the reflection principle to
Brownian bridge; see also \cite[Lemma 9]{bramson2}.
  (The bound in discrete time is
the same as that in continuous time, up to the constant $C_{25}$, since the
``overshoot" of the normal past a boundary, over a unit time
interval, has bounded
second moment.)

%Bound in Proposition 6.3.  
\subsubsection*{Bound in Lemma \ref{lem-11}}
The bound in (6.9) is obtained in the same manner as
are parts (a) and (b) of Lemma 11 on page 565 of \cite{bramson2}.  
One can apply the
reflection principle as in part (a), but for discrete time instead of
continuous time.  (As in the derivation of the upper bound of 
(\ref{eq-260810c}), the
``overshoot" of the normal only affects the constant in front.)
One obtains with a little work that, for
$w\le r(z,z^{\prime})(A_{n}+1)/n - L_{n}(z,z^{\prime})+1$,
\begin{equation}
\label{eqA3}
P(D_{z,z^{\prime},w}^{(2)}) \le C_{26} n^{-\frac{3}{2}}(-w +
\frac{r(z,z^{\prime})}{n}A_n + 2) \exp\{-(A_n - w)^2 / 2u(z,z^{\prime})\}.
\end{equation}

As in part (b) of Lemma 11 of \cite{bramson2}, 
for the above range of $w$, the right side of (\ref{eqA3})
is maximized at the boundary $w = r(z,z^{\prime})(A_{n}+1)/n -
L_{n}(z,z^{\prime})+1$.  Plugging this value of $w$ into the right side
of
(\ref{eqA3}) yields (\ref{eq-270810f}).

\end{appendix}

\begin{tabular}{llll}  & Maury Bramson&& Ofer Zeitouni\\
  & School of
Mathematics &&  School of Mathematics\\
&University of Minnesota&& University of Minnesota\\
&206 Church St. SE && 206 Church St. SE\\
&Minneapolis, MN 55455&& Minneapolis, MN 55455\\
&bramson@math.umn.edu&& \; and \\  
& 
&& Faculty of Mathematics\\
 &&& Weizmann Institute of Science\\
 &&& Rehovot 76100, Israel\\
&&& zeitouni@math.umn.edu
\end{tabular}

\end{document}